\documentclass[10pt]{amsart}
\usepackage{eucal, mathrsfs}
\usepackage[OT2,OT1,T1]{fontenc}
\usepackage{amsmath,amssymb,amsthm,amsfonts,amscd,amsopn}
\usepackage{hyperref}
\usepackage{newcent}
\usepackage{epsfig}
\usepackage{xspace}
\usepackage[all]{xy}
\xyoption{all}
\def\cyr{\fontencoding{OT2}\fontfamily{wncyr}\selectfont}

\usepackage{fullpage}

\newcommand{\field}[1]{\mathbf #1}
\newcommand{\mf}[1]{\mathfrak #1}

\newcommand{\ms}[1]{\mathscr #1}
\newcommand{\widebar}[1]{\overline{#1}}

\newcommand{\R}{\field R}

\newcommand{\C}{\field C}
\newcommand{\F}{\field F}
\newcommand{\Z}{\field Z}
\newcommand{\Q}{\field Q}

\newcommand{\simto}{\stackrel{\sim}{\to}}

\renewcommand{\phi}{\varphi}

\renewcommand{\hom}{\operatorname{Hom}}
\newcommand{\uhom}{\operatorname{\underline{Hom}}}
\newcommand{\shom}{\ms H\!om}
\DeclareMathOperator{\chom}{\underline{Hom}}

\newcommand{\saut}{\ms A\!ut}

\newcommand{\send}{\ms E\!nd}

\newcommand{\spec}{\operatorname{Spec}}

\DeclareMathOperator{\Pic}{Pic}
\newcommand{\sPic}{\ms P\!ic}
\newcommand{\td}{\operatorname{Td}}
\DeclareMathOperator{\chern}{ch}

\newcommand{\m}{\boldsymbol{\mu}}

\newcommand{\G}{\field G} 

\renewcommand{\H}{\operatorname{H}}

\newcommand{\GL}{\operatorname{GL}}

\DeclareMathOperator{\Ext}{\operatorname{Ext}}

\newcommand{\Sha}{\text{\cyr Sh}}

\DeclareMathOperator{\Coh}{\operatorname{Coh}}

\newcommand{\Gal}{\operatorname{Gal}}

\DeclareMathOperator{\per}{per}
\DeclareMathOperator{\ind}{ind}

\DeclareMathOperator{\D}{\operatorname{\bf D}}

\DeclareMathOperator*{\tensor}{\otimes}

\DeclareMathOperator{\rk}{\operatorname{rk}}

\newcommand{\surj}{\twoheadrightarrow}

\newcommand{\id}{\operatorname{id}}

\DeclareMathOperator{\End}{\operatorname{End}}
\DeclareMathOperator{\Aut}{\operatorname{Aut}}
\DeclareMathOperator{\aut}{\operatorname{Aut}}
\DeclareMathOperator{\isom}{\operatorname{Isom}}
\DeclareMathOperator{\Isom}{\operatorname{Isom}}
\DeclareMathOperator{\M}{\operatorname{M}}
\DeclareMathOperator{\Br}{\operatorname{Br}}

\DeclareMathOperator{\B}{\operatorname{B\!}}

\newcommand{\HH}{\mathbf H}
\newcommand{\FF}{\mathbf F}

\DeclareMathOperator{\cSimp}{\underline{Spl}}
\DeclareMathOperator{\Simp}{Spl}

\newtheorem{lem}{Lemma}[subsection]
\newtheorem{thm}[lem]{Theorem}

\newtheorem{prop}[lem]{Proposition}

\newtheorem{cor}[lem]{Corollary}

\theoremstyle{definition}
\newtheorem{defn}[lem]{Definition}
\newtheorem{con}[lem]{Convention}
\newtheorem{example}[lem]{Example}

\theoremstyle{remark}
\newtheorem{remark}[lem]{Remark}

\newtheorem{ques}[lem]{Question}

\DeclareMathOperator{\cpc}{cap}

\newcommand{\modify}[1]{\widetilde{#1}}

\newtheorem{btheorem}{Theorem}[subsection]
\newtheorem{blemma}[btheorem]{Lemma}
\newcommand{\NS}{{\mathrm{NS}}}

\title{Index reduction for Brauer classes via stable sheaves}
\author{Daniel Krashen}
\address{Yale University, New Haven, CT}
\email{daniel.krashen@yale.edu}
\thanks{Krashen was partially supported by an NSA Young Investigator's Grant.}
\author{Max Lieblich}
\address{Princeton University, Princeton, NJ}
\email{lieblich@math.princeton.edu}
\thanks{Lieblich was partially supported by an NSF Postdoctoral Fellowship.}
\address{Princeton University, Princeton, NJ}
\email{bbhatt@math.princeton.edu}

\begin{document}

\begin{abstract}
  We use twisted sheaves to study the problem of index reduction for
  Brauer classes.  In general terms, this problem may be phrased as
  follows: given a field $k$, a $k$-variety $X$, and a class $\alpha
  \in \Br(k)$, compute the index of the class $\alpha_{k(X)} \in
  \Br(X)$ obtained from $\alpha$ by extension of scalars to $k(X)$. We
  give a general method for computing index reduction which refines
  classical results of Schofield and van den Bergh.  When $X$ is a
  curve of genus $1$, we use Atiyah's theorem on the structure of
  stable vector bundles with integral slope to show that our formula
  simplifies dramatically, giving a complete solution to the index
  reduction problem in this case.  Using the twisted Fourier-Mukai
  transform, we show that a similarly simple formula describes
  homogeneous index reduction on torsors under higher-dimensional
  abelian varieties.
\end{abstract}
\maketitle
\tableofcontents
\section{Introduction}
\label{sec:introduction}

Let $k$ be a field and $X$ a $k$-variety. The index reduction problem
asks the following question: Given an element of the Brauer group
$\alpha \in \Br(k)$, how does one compute the index of the class
$\alpha_{k(X)}$ obtained by extending scalars from $k$ to $k(X)$? The
implicit assumption in this question is that one may already know how to
compute the index of an algebra defined over $k$ (and its finite field
extensions). A reasonable answer therefore would be to describe the
index of $\alpha_{k(X)}$ in terms of the indicies of some list of
auxiliary algebras defined over the ground field $k$, or finite
extensions of $k$.  Results of this sort have been obtained for certain
hyperelliptic curves $X$ and $k$ a p-adic field in \cite{Yanchevskii},
and for $k$ arbitrary and $X$ a projective homogeneous variety for some
linear algebraic group in \cite{Merk-Panin-Wadsworth}. In this paper we
give a solution in the case that $X$ is a curve of genus $1$ (Corollary
\ref{sec:index-reduction-via-3}), and make some reductions on the
computation of the index in the general case (Proposition
\ref{sec:index-reduction-via-2}).

Our method is based on the use of twisted
sheaves. Before going into detail about how the theory is used in this
particular case, it may be instructive to give a philosophical
comparison between our method and the methods used in
\cite{Merk-Panin-Wadsworth} in the case of homogeneous varieties.

To begin, we consider a Brauer class on a variety which has been
obtained by extension of scalars from the ground field. In the above
notation, this would be the class $\alpha_X$. Finding the index of
$\alpha_{k(X)}$ corresponds to the finding the minimal dimension of a
module for an Azumaya algebra representing the class $\alpha_{k(X)}$.
Concretely, if $\alpha_{k(X)}$ is the class of a central simple algebra
$A = M_m(D)$ for some $k(X)$-central division algebra $D$, then a module
of minimal dimension would be of the form $D^m$ and we may compute the
index of $A$ as $\ind(A) = \deg(D) = \dim_k(D^m)/\deg(A)$.  The problem
of finding a module of minimal size may be made more geometric by
noticing that this module corresponds to a coherent sheaf of modules for
an Azumaya algebra $A_X$ in the class $\alpha_X$. In
\cite{Merk-Panin-Wadsworth} the authors then proceed by considering the
$K$-theory of the category of such modules.

The main idea in our approach here is based on the desire to deal with
sheaves of modules for $A_X$ is a more geometric way, and in particular,
in a way which reflects more closely the study of sheaves of modules
over $X$, allowing us to draw from the theory of vector bundles on $X$.
To accomplish this, the algebra $A_X$ is replaced by a gerbe
corresponding to its Brauer class $\alpha_X$, and correspondingly the
sheaves of modules for $A_X$ are replaced by twisted sheaves for the
gerbe. This allows us to realize our goal, as many useful facts about
vector bundles and sheaves turn out to generalize to their twisted
counterparts. In the case of genus $1$ curves, Atiyah's classification
of vector bundles on an elliptic curve plays an essential role in the
problem of index reduction.

Attempting to generalize these results to torsors under
higher-dimensional abelian varieties leads us to briefly study the
twisted Fourier-Mukai transform.  We deduce a criterion for
homogeneous index reduction (see \ref{D:homog-ind-red}) which shows
the stark difference between index reduction for torsors under abelian
varieties and rational homogeneous spaces. We include an appendix by
Bhargav Bhatt which uses the twisted Fourier-Mukai transform to study
the period-index problem for torsors under abelian varieties.

In this paper we freely use the theory of twisted sheaves.  For the
basic facts concerning these objects and their relation to the index of
a Brauer class, the reader should refer to \cite{period-index-paper}.

\section*{Acknowledgments}

During the course of this work, we had helpful conversations with
Bhargav Bhatt and Jean-Louis Colliot-Th\'el\`ene.

\subsection{Simple bundles on a pullback gerbe}
\label{sec:simple-bundles}

In laying the groundwork for our constructions, it is not necessary to
work over a field, and so we develop our basic machinery in the context
of a base $S$ which is an algebraic space. If the reader prefers, they
may simply consider the case $S = \spec(k)$.

Let $f:X\to S$ be a proper flat morphism of finite presentation
between algebraic spaces which is cohomologically flat in dimension
$0$ and $\ms X\to X$ a $\m_n$-gerbe.  We begin by introducing the
stack $\cSimp_{\ms X/S}$ of simple $\ms X$-twisted sheaves and its
corresponding coarse moduli space $\Simp_{\ms X/S}$.  As we will see
below (~\ref{sec:simple-bundl-pullb-1}), the natural map $\cSimp_{\ms
  X/S} \to \Simp_{\ms X/S}$ is a $\G_m$-gerbe, and thus thus the
obstruction for lifting a point in the coarse moduli to an object of
the stack may be interpreted as a Brauer class.

\begin{defn}\label{sec:simple-bundl-pullb-2}
  Given $T\to S$, a \emph{flat family of simple coherent $\ms
    X$-twisted sheaves parametrized by $T$\/} is a $T$-flat
  quasi-coherent $\ms X$-twisted sheaf of finite presentation $\ms F$
  on $X_T$ such that
  \begin{enumerate}
  \item the formation of $(f_T)_{\ast}\send(\ms F)$ commutes with base
    change on $T$, and
  \item the natural map $\ms O_T\to (f_T)_{\ast}\send(\ms F)$ is an
    isomorphism.
  \end{enumerate}
\end{defn}

\begin{lem}\label{basechange}
  A flat family of coherent $\ms X$-twisted sheaves $\ms F$
  parametrized by $T$ is simple if and only if for all geometric
  points $\widebar t\to T$ the fiber $\ms F_{\widebar t}$ is simple.
\end{lem}
\begin{proof}
  If the family is simple then it is simple on each fiber by
  definition.  Suppose that each geometric fiber of $\ms F$ is simple.
  To prove that $\ms F$ is a simple family, it suffices to show that
  the natural map $\sigma:\ms O_T\to(f_T)_{\ast}\send(\ms F)$ is an
  isomorphism under the assumption that $T$ is a local Noetherian
  scheme and the geometric closed fiber $\ms F_{\widebar t}$ is
  simple.  Since $\ms F_{\widebar t}$ is simple, $\ms F$ has
  non-trivial restriction to the closed fiber, we see that $\ms F$ is
  faithfully flat over $T$.  This implies that $\sigma$ is injective,
  so it remains to show that $\sigma$ is surjective.  To show this, we
  may assume (by the faithful flatness of completion) that $T$ is the
  spectrum of a complete local Noetherian ring $A$ with maximal ideal
  $\mf m$.  The Grothendieck Existence Theorem implies that the
  natural map
$$\End(\ms F)\simto\lim \End(\ms F\tensor A/\mf m^n\ms F)$$
is an isomorphism of $A$-modules, and this is clearly compatible with
the same natural isomorphism
$$A\simto\lim A/\mf m^n.$$
Thus, it suffices to show the statement assuming that $T$ is an
Artinian local ring; let $k$ be the residue field of $A$ and $\mf m$
the maximal ideal.  Letting $\mf m^n=0$ be the minimal power of the
maximal ideal which vanishes, we may assume by induction (and the
assumption on the fiber) that $\sigma_{n-1}$ is an isomorphism, where
$\sigma_{n-1}:A/\mf m^{n-1}\to\End(\ms F\tensor A/\mf m^{n-1})$ is the
natural map.  Standard results in deformation theory show that the map
$f\mapsto f\tensor_k\mf m^{n-1}$ gives the kernel of the natural
restriction map, yielding an exact sequence
$$0\to \End(\ms F\tensor k)\tensor\mf m^{n-1}\to\End(\ms F)\to\End(\ms F\tensor A/\mf m^{n-1}).$$
This sequence admits a map from the exact sequence
$$0\to\mf m^{n-1}\to A\to A/\mf m^{n-1}\to 0$$
which, by induction, is an isomorphism on the outer non-zero terms.
It follows that the central term is an isomorphism, as desired.
\end{proof}

It is clear that the collection of flat families of simple coherent
$\ms X$-twisted sheaves is a stack in the fpqc topology on
$S$-schemes.  We will write $\cSimp_{\ms X/S}$ for this stack.

\begin{lem}\label{sec:simple-bundl-pullb}
  The inertia stack $\ms I(\cSimp_{\ms X/S})$ is naturally isomorphic
  to $\G_{m,\cSimp_{\ms X/S}}$.
\end{lem}
\begin{proof}
  With the notation of \ref{sec:simple-bundl-pullb-2}, the map $\ms
  O_T\to(f_T)_{\ast}\send(\ms F)$ gives rise to a map
  $\G_{m,T}\to(f_T)_\ast\Aut(\ms F)$.  This yields a natural map $\G_m\to\ms
  I(\cSimp_{\ms X/S})$.  The conditions of
  \ref{sec:simple-bundl-pullb-2} then say precisely that this map is
  an isomorphism.
\end{proof}

\begin{lem}\label{sec:simple-bundl-pullb-3}
  Let $\ms F$ be a $T$-flat quasi-coherent sheaf of finite
  presentation on $X_T$.  There is an open subscheme $U\subset T$
  such that a map $T'\to T$ factors through $U$ if and only if the
  pullback $\ms F_{T'}$ is a flat family of simple coherent $\ms
  X$-twisted sheaves.
\end{lem}
\begin{proof}
  Since being simple is a fiberwise condition, it suffices to assume
  that $T$ is reduced and Noetherian and show that the set $U$ of points
  parametrizing simple fibers is open.  The set is constructible:
  $f_{\ast}\send(\ms F)$ is coherent and generically compatible with
  base change, as are the kernel and cokernel of the map $\ms O_T\to
  f_{\ast}\send(\ms F)$.  Nakayama's lemma immediately shows
  constructibility.  The set is stable under generization: suppose $T$
  is the spectrum of a discrete valuation ring whose closed point is in
  $U$.  The formation of $f_{\ast}\send(\ms F)$ is certainly compatible
  with passage to the generic fiber.  On the other hand,
  $f_{\ast}\send(\ms F)$ is a coherent sheaf on $T$ whose closed fiber
  is $1$-dimensional.  It follows from standard semicontinuity results
  (for possibly non-flat modules over a dvr) that the generic fiber is
  also $1$-dimensional.  From this it immediately follows that the
  scalars generate the endomorphisms on the generic fiber.  Since $U$ is
  constructible and stable under generalization, we conclude that it is
  open, as desired.
\end{proof}

\begin{prop}\label{sec:simple-bundl-pullb-1}
  The stack $\cSimp_{\ms X/S}$ is a $\G_m$-gerbe over an algebraic
  space locally of finite presentation $\Simp_{\ms X/S}\to S$.
\end{prop}
\begin{proof}
   The methods of \cite{twisted-moduli} show that the stack $\Coh_{\ms
     X/S}$ of flat families of coherent sheaves on $\ms X$ is an Artin stack locally
   of finite presentation over the base.  It is straightforward to
   check that the collection of $\ms X$-twisted coherent sheaves is an
   open substack.  Applying \ref{sec:simple-bundl-pullb-3}, we see that
   $\cSimp_{\ms X/S}$ is an open substack of $\Coh_{\ms X/S}$, and thus
   it is an Artin stack locally of finite presentation over $S$.
 
   Since the inertia stack is fppf over $\cSimp_{\ms X/S}$, it follows
   from standard methods (first described in the last paragraph of
   remark 2 in the appendix to \cite{artin}) that the sheafification
   $\Simp_{\ms X/S}$ of $\cSimp_{\ms X/S}$ is an algebraic space and
   that the natural map $\cSimp_{\ms X/S}\to\Simp_{\ms X/S}$ realizes
   $\cSimp_{\ms X/S}$ as a gerbe.  It is then immediate that it is a
   $\G_m$-gerbe.
\end{proof}

\begin{defn}\label{sec:simple-bundl-pullb-5}
  With the above notation, the cohomological Brauer class
  $[\cSimp_{\ms X/S}]\in\H^2(\Simp_{\ms X/S},\G_m)$ will be called the
  \emph{universal obstruction\/}.
\end{defn}

When the gerbe $\ms X\to X$ is a pullback from $S$, we can be more
precise about the structure of the universal obstruction and its
variation with $\ms S$.

\begin{prop}\label{sec:simple-bundl-pullb-4}
  Let $\ms S\to S$ be a $\m_n$-gerbe, and let $\modify{\ms S} \to S$ be
  the associated $\G_m$-gerbe.
  \begin{enumerate}
  \item There is a natural isomorphism
    $\phi_{\ms S}:\Simp_{X/S}\simto\Simp_{X\times_S\ms S}$ over $S$.
  \item Via $\phi$, there is an equality $$[\cSimp_{X\times_S\ms
      S}]-[\cSimp_{X/S}]=[\modify{\ms S}_{\Simp_{X/S}}]\in\H^2(\Simp_{X/S},\G_m).$$
  \end{enumerate}
\end{prop}
\begin{proof}
  We will define a section of the sheaf
  $\Isom_S(\Simp_{X/S},\Simp_{X\times_S\ms S/S})$ by gluing local
  sections.  If $\modify{\ms S}$ is trivial, so that there is an invertible
  $\ms S$-twisted sheaf $\ms L$, then $V\mapsto V\tensor\ms L$ defines
  a diagram $$\xymatrix{\cSimp_{X/S}\ar[r]\ar[d] &
    \cSimp_{X\times_S\ms S/S}\ar[d]\\ \Simp_{X/S}\ar[r]^-{\phi} &
    \Simp_{X\times_S\ms S/S}.}$$ Thus, on an \'etale cover $U\to S$
  such that $\modify{\ms S}\times_S U$ is trivial, we have such a section
  $\phi_U$.  On $U\times_S U$ we have $p_1^{\ast}\ms L\tensor
  p_2^{\ast}\ms L^{\vee}\cong\ms M$ with $\ms M\in\Pic(U\times_S U)$.
  Thus, $p_2^{\ast}\phi_U^{-1}\circ p_1^{\ast}\phi_U$ is the map
  $\Simp_{X/S}\to\Simp_{X/S}$ induced by twisting by $\ms M$.  On the
  coarse moduli space, this map equals $\id$.  Thus, $\phi_U$
  satisfies the cocycle condition, giving rise to a global
  isomorphism.  A similar argument shows that the resulting
  isomorphism is independent of the choices of $U$ and $\ms L$.

  To prove the second statement, we use the results of section 2.4 of
  chapter IV of \cite{giraud}.  First, we note that the stack
  $\modify{\ms S}$ parametrizes $\ms S$-twisted invertible sheaves.
  (In other words, the stack $\chom_S^{\G_m}(\modify{\ms S},\B{\G_m})$
  of 1-morphisms of $\G_m$-gerbes is isomorphic to $\modify{\ms S}$.)
  Consider the morphism
  $$\cSimp_{X/S}\times_S \ms S \cong 
  \cSimp_{X/S}\times_{\Simp_{X/S}}\ms S_{\Simp_{X/S}} \to
  \cSimp_{X/S}\times_{\Simp_{X/S}} \modify{\ms S}_{\Simp_{X/S}} \to
  \cSimp_{X\times_S \modify{\ms S} /S}$$
  which sends a pair $(V,L)$ consisting of a flat family of simple
  coherent sheaves on $X$ and a flat family of invertible $\ms
  S$-twisted sheaves to the tensor product $V\tensor L$.  This map is
  compatible with the product map $\G_m\times\G_m\to\G_m$ in the obvious
  way.  In the notation of section 1.6.1 of chapter IV, the
  multiplication map is precisely the contracted product of the two
  copies of $\G_m$.  By Proposition IV.2.4.1 of \cite{giraud}, there is an
  induced map of $\G_m$-gerbes
  $$\cSimp_{X/S}\bigwedge^{C}_{\Simp_{X/S}} \modify{\ms S}_{\Simp_{X/S}}
  \to\cSimp_{X\times_S \modify{\ms S}/S}.$$
  But any map of gerbes gives an equality of the associated cohomology
  classes.  Since the contracted product gives the sum of cohomology
  classes, the result follows.
\end{proof}

\begin{cor}\label{C:tw-pic-obs-ch}
  The twisted Picard stack $\sPic_{X\times_S\ms S/S}$
  naturally has sheafification
  $\Pic_{X/S}$ with universal obstruction $[\sPic_{X/S}]+[\modify{\ms S}_{\Pic_{X/S}}]$.
\end{cor}

\subsection{(Semi)stable bundles}
\label{sec:semistable-bundles}

We briefly recall the elements of the theory of stable and semistable
twisted sheaves on curves.  This theory is the specialization of a
much more general theory of (semi)stable sheaves on arbitrary
polarized gerbes (or even arbitrary polarized orbifolds), but such
extreme generality will have no place in the rest of this paper.  We
refer the reader to the first chapters of \cite{h-l} for more on the
general theory.

Let $X/k$ be a proper smooth geometrically connected curve over a
field and let $\pi:\ms X\to X$ be a $\m_n$-gerbe.  If $L$ is an
invertible sheaf on $\ms X$, then the natural map
$\pi^{\ast}\pi_{\ast}(L^{\tensor n})\to L^{\tensor n}$ is an
isomorphism.  This permits us to define a degree for invertible
sheaves on the stack $\ms X$ (and by extension, if necessary, the
degree of any coherent sheaf.)

\begin{defn}
  Given an invertible sheaf $L$ on $\ms X$, the \emph{degree\/} of $L$
  is $$\deg(L)=\frac{1}{n}\deg_X(\pi_{\ast}(\ms L^{\tensor
    n}))\in\Q.$$  Given a locally free sheaf $V$ on $\ms X$, the
  \emph{slope\/} of $V$ is $\mu(V)=\deg(\det V)/\rk V$.
\end{defn}
Using the fact that every coherent sheaf on
$\ms X$ has a finite resolution by locally free sheaves, one can
extend the definition of degree and slope to arbitrary coherent
sheaves.  While this is entirely reasonable, it will not come up in
the sequel.

\begin{defn}\label{sec:semistable-bundles-3}
  A sheaf $V$ on $\ms X$ is \emph{stable\/} (resp.\
  \emph{semistable\/}) if it is locally free and for all proper
  subsheaves $F\subsetneq V$ one has $\mu(F)<\mu(V)$ (resp.\
  $\mu(F)\leq\mu(V)$).
\end{defn}

\begin{remark}\label{sec:semistable-bundles-2}
  It is easy to check that the condition that $\mu(F)\leq\mu(V)$ is
  the same as $\mu(V)\leq\mu(V/F)$.  We will use this implicitly in
  the sequel.
\end{remark}

One can similarly define the notion of \emph{geometrically
(semi)stable\/}.  It turns out that semistability is a geometric
property, while stability is not unless the sheaf in question is
simple (see Example 1.3.9 and \S{1.5} of \cite{h-l}).  It is easy
to see that one can test the property of (semi)stability by
restricting attention solely to subsheaves $F\subset V$ for which the
quotient $V/F$ is also locally free.

\begin{defn}
  A semistable sheaf $V$ on $\ms X$ is \emph{(geometrically)
    polystable\/} if $V$ (resp.\ $V\tensor\widebar k$) is isomorphic
  to a direct sum of stable sheaves.
\end{defn}

If $V$ is polystable, it follows from Remark
\ref{sec:semistable-bundles-2} that the stable summands all have the
same slope, and that this slope equals the slope of $V$.

\begin{remark}
  When the gerbe $\ms X$ is trivial, it is in general quite subtle to
  detect semistable sheaves.  However, if $\ms X$ is non-trivial, then
  any locally free $\ms X$-twisted sheaf $V$ whose rank equals the
  index of the Brauer class attached to $\ms X$ is automatically
  stable.  However, in this case it is quite subtle to detect when
  such a sheaf is geometrically stable.  This issue will appear in
  a fundamental way in the analysis of section \ref{sec:index-reduction-via}.
\end{remark}

\begin{prop}\label{sec:semistable-bundles-1}
  Given $\mu\in\Q$, the category of semistable $\ms
  X$-twisted sheaves of slope $\mu$ is an abelian category in which every
  object has finite length.
  The simple objects are the stable sheaves. 
\end{prop}
\begin{proof}
  Let $\phi:F\to G$ be a map of semistable sheaves of slope $\mu$.  We
  will show that the kernel $K$ and cokernel $Q$ of $\phi$ are both
  semistable of slope $\mu$.  This will immediately prove the first
  statement by ``transport of structure'' from the abelian category of
  coherent sheaves on $\ms X$.  We first claim that $Q$ is locally free.
  Indeed, we have 
$$\mu=\mu(F)\leq\mu(\phi(F))\leq\mu(\phi(F)^\ast)\leq\mu(G)=\mu,$$
where $\phi(F)^\ast$ denotes the saturation of $\phi(F)$ as a subsheaf
of $G$.  We conclude that $\phi(F)=\phi(F)^\ast$, so that
$G/\phi(F)\cong G/\phi(F)^\ast$.  It follows that $Q$ is locally free.
Moreover, we have that the first two non-trivial terms in the sequence
$0\to\phi(F)\to G\to Q\to 0$ have the same slope.  It follows that
$\mu(Q)=\mu$.  To see that $Q$ is semistable, suppose $Q'\subset Q$ is
a subsheaf and let $G'$ be the preimage of $Q'$ in $G$, so that there
is an exact sequence $0\to\phi(F)\to G'\to Q'\to 0$.  By the
semistability of $G$ we have that $\mu(\phi(F))\geq\mu(G')$, which
implies that $\mu(G')\geq\mu(Q')$.  Thus, $\mu(Q')\leq\mu$ and $Q$ is
semistable.

It is immediate that $K$ is locally free.  Since
$\mu(F)=\mu(\phi(F))$, we conclude that $\mu(K)=\mu(F)=\mu$.  It now
follows immediately from the semistable of $F$ that $K$ is semistable.
\end{proof}

\begin{cor}\label{sec:semistable-bundles-4}
  If $\ms F$ is semistable of slope $\mu$ and $\{\ms G_i\}_{i\in I}$ is a
  set of stable subsheaves of $\ms F$ of slope $\mu$ then there is a
  subset $J\subset I$ such that $\Sigma_{i\in I}\ms
  G_i=\bigoplus_{j\in J}\ms G_j$ as subsheaves of $\ms F$.  
\end{cor}
\begin{proof}
  Since $\ms F$ is Noetherian, we may assume that $I$ is finite.
  Consider the surjection $\bigoplus_i\ms G_i\surj\Sigma_i\ms G_i$,
  and write $K$ for the kernel.  Let $K'\subset K$ be a stable
  subsheaf of slope $\mu$.  The map $K'\to\bigoplus\ms G_i$ is
  non-zero, which implies that one of the projections $K'\to\ms
  G_{i_0}$ is an isomorphism.  It follows that $\Sigma_i\ms G_i$ is
  the image of $\bigoplus_{i\neq i_0}\ms G_i$.  By induction,
  $\Sigma_i\ms G_i$ is a direct sum of stable subsheaves.
\end{proof}

The \emph{socle\/} of a semistable sheaf $\ms F$ of slope $\mu$ is the sum of
all of its stable subsheaves of slope $\mu$.  It is clear that the
socle is stable under all automorphisms of $\ms F$, and we have just
shown that the socle is polystable.

\begin{cor}
  Let $k$ be a perfect field. Given a semisimple $\ms X$-twisted sheaf
  $\ms F$ of slope $\mu$, there is a canonical maximal subsheaf $S(\ms
  F)\subset\ms F$, compatible with base extension, such that $S(\ms F)$
  is geometrically polystable.
\end{cor}
\begin{proof}
  It follows from \ref{sec:semistable-bundles-1} that the sum 
  $\Sigma\ms G\subset\ms F\tensor\widebar k$, taken over all subsheaves
  $\ms G\subset\ms F\tensor k$ which are stable of slope $\mu$, is
  direct and stable under the operation of the absolute Galois group
  of $k$.  Basic descent theory shows that this sheaf is the base
  extension of a subsheaf $S(\ms F)\subset\ms F$.  Uniqueness is immediate.
\end{proof}

\subsection{Moduli of stable sheaves}
\label{sec:moduli-stable-sheav-1}

Throughout this section, $X$ will be a smooth proper geometrically
curve over a \textit{perfect} field $k$.  We recall from section 2.3
of \cite{twisted-moduli} that the $\ms X$-twisted stable sheaves on
$X$ of slope $\mu$ form an algebraic Deligne-Mumford stack which we
denote by $\ms M^{s, \mu}_{\ms X/k}$.  This is naturally a substack of
$\cSimp_{\ms X/k}$. We denote its coarse moduli space by $\M^{s,
  \mu}_{\ms X/k}$ and note that $\ms M^{s, \mu}_{\ms X/k}$ is a
$\G_m$-gerbe and may in fact be obtained as the pullback of the morphism 
$\cSimp_{\ms X/k} \to \Simp_{\ms X/k}$ via the inclusion $\M^{s,
  \mu}_{\ms X/k} \to \Simp_{\ms X/k}$.

We will now study the possible ranks of certain geometrically
poly\-stable sheaves. Given a geometrically polystable $\ms X$-twisted
sheaf $V$ on $X$ of slope $\mu$, we have a decomposition
$V\tensor\widebar k\cong\bigoplus V_i$ with $V_i$ stable sheaves of
slope $\mu$ on $X\tensor\widebar k$.  Each $V_i$ gives rise to a point
$[V_i]$ on the coarse moduli space $\M^{s,\mu}_{\ms X/k}$ of stable $\ms
X$-twisted sheaves of slope $\mu$.  Write $I_V$ for the set of points
$[V_i]\in\M^{s,\mu}_{\ms X/k}(\widebar k)$.

\begin{lem}\label{sec:moduli-stable-sheav-3}
  With the above notation, there is a natural continuous action of
  $\Gal(k)$ on $I_V$.
\end{lem}
\begin{proof}
  Write $V\tensor\widebar k=\bigoplus_{i\in I_V}W_i$, where $W_i$
  groups the stable summands with the given isomorphism class.  Given
  an element $\sigma\in\Gal(k)$, the descent datum on
  $V\tensor\widebar k$ induces an isomorphism
  $\bigoplus\sigma^{\ast}W_i\simto\bigoplus W_i$.  Since each $W_i$
  (resp.\ $\sigma^{\ast}W_i$) is an isotypic and the $W_i$ (resp.\
  $\sigma^{\ast}W_i$) have pairwise non-isomorphic stable
  constituents, it follows that there is an induced bijection
  $\widetilde{\sigma}:I_V\simto I_V$ such that the isomorphism class
  associated to $\sigma^{\ast}W_i$ is the same as that associated to
  $W_{\widetilde{\sigma}(i)}$.  The map $\sigma\mapsto\widetilde{\sigma}$
  defines the action in question.
\end{proof}

\begin{lem}\label{sec:moduli-stable-sheav-4}
  In the notation of \ref{sec:moduli-stable-sheav-3}, the sheaf $V$ is
  indecomposable if and only if the action of $\Gal(k)$ on $I_V$ is transitive.
\end{lem}
\begin{proof}
  The action of $\Gal(k)$ on an orbit of $I_V$ induces (via
  restriction) a descent datum on a proper sub-sum of the $W_i$.
  Taking the direct sum over all orbits yields a decomposition of $V$
  as a direct sum indexed by orbits.
\end{proof}

\begin{lem}\label{sec:moduli-stable-sheav-6}
  Let $k$ be a field and $\alpha$ and $\beta$ two elements of
  $\Br(k)$.  Suppose that for all field extensions $L/k$, we have that
  $\alpha_L=0$ if and only if $\beta_L=0$.  Then
  $\alpha$ and $\beta$ generate the same cyclic subgroup of $\Br(k)$.
\end{lem}
\begin{proof}[Sketch of proof]
  This is a well-known result of Amitsur.  We provide a modern proof.
  Let $P$ be a Brauer-Severi variety with Brauer class $\alpha$.
  Since $\alpha|_P=0$, we have that $\beta|_P=0$.  Examining the Leray
  spectral sequence in \'etale cohomology for $\G_m$ on the morphism
  $P\to\spec k$, we see that the kernel of $\Br(k)\to\Br(P)$ is the
  subgroup generated by $\alpha$.  Thus,
  $\beta\in\langle\alpha\rangle$.  Reversing the roles of $\alpha$ and
  $\beta$ completes the proof.
\end{proof}

\begin{prop}\label{sec:moduli-stable-sheav-5}
  Given an indecomposable geometrically polystable $\ms X$-twisted
  sheaf $V$, there is a $\mu \in \Q$ and a closed point $p\in\M^{s,
  \mu}_{\ms X/k}$ such that the rank of $V$ is a multiple of
  $[\kappa(p):k]\ind(\alpha(p))$.  Moreover, every such multiple is
  realized by an indecomposable geometrically polystable $\ms X$-twisted
  sheaf.
\end{prop}
\begin{proof}
  Write $V\tensor\widebar k=\bigoplus V_i=\bigoplus W_j$ as above.  By
  \ref{sec:moduli-stable-sheav-4}, the points $[V_i]$ form a single
  Galois orbit in $\M^{s, \mu}_{\ms X/k}(\widebar k)$.  This corresponds to a
  closed point $p\in\M^{s, \mu}_{\ms X/k}$.

  Consider the $k$-algebra $A:=\End(V)$.  Tensoring with $\widebar k$,
  we see that $A$ is a separable algebra.  Since $V$ is indecomposable,
  it is easy to see that $A$ is in fact a division algebra.  Thus, there
  is a finite extension $L$ of $k$ such that $A$ is a central division
  algebra over $L$.  We claim that $L=\kappa(p)$ and that the class of
  $A$ in $\Br(L)$ equals the universal obstruction $\omega(p)$
  restricted to $p$.

  To prove that $L=\kappa(p)$, note that the set of idempotents of
  $A\tensor\widebar k$ is in natural bijection with $I_V$, in a manner
  compatible with the action of $\Gal(k)$.  It follows that the finite
  \'etale coverings $\spec\kappa(p)\to\spec k$ and $\spec Z(A)\to\spec
  k$ are isomorphic, which shows that $L\cong\kappa(p)$.  To see that
  $[A]=\omega(p)$, it suffices to do so assuming that $k=L$.  Indeed,
  the inclusion $L\subset\End(A)$ gives $V$ the structure of $X\tensor
  L$-module in such a way that $V\tensor_L\widebar L$ is polystable
  and isotypic.  Thus, we may assume that $A$ is a central division
  algebra over the base field and that $V\tensor_k\widebar k$ is
  isotypic.  

  To show that $[A]=\omega(p)$, we first show that $[A]$ and
  $\omega(p)$ generate the same cyclic subgroup of $\Br(k)$. By
  \ref{sec:moduli-stable-sheav-6}, it suffices to show (upon extending
  the base field, which we will denote with $k$ by abuse of notation)
  that $A$ is split if and only if there is a stable bundle $V_0$ on
  $X$ with moduli point $p$.  If $V_0$ exists, then we see that $V$ is
  an \'etale form of $V_0^{\oplus n}$ for some $n$.  But
  $\Aut(V_0^{\oplus n})=\GL_{n,k}$, so by Hilbert's Theorem 90 $V\cong
  V_0^{\oplus n}$.  Thus, $A\cong\M_n(k)$.  On the other hand, if $A$
  is split then there is a full set of operators on $V$ splitting $V$
  into a direct sum $V_0^{\oplus n}$ with $V_0$ geometrically stable.
  It follows that $[V_0]=p$ and thus that $\omega(p)=0$.  (The reader
  will note that for the purposes of this paper, the equality of the
  cyclic subgroups $\langle[A]\rangle$ and $\langle\omega(p)\rangle$
  is all that we need.)

  Note that given any $V$, its forms are classified by
  $\H^1(\spec k,\GL(A))$.  By Hilbert's Theorem 90 for division
  algebras, this cohomology group is trivial.  Thus, any two
  geometrically polystable geometrically isotypic $\ms X$-twisted
  sheaves with geometric summands supported at $p$ are isomorphic.
  Furthermore, the argument of the previous paragraph shows that the
  rank of $V$ is a multiple of $r\ind(\omega(p))$.  Hence, to show
  that $[A]=\omega(p)$ it suffices to show that there is a $V$ of rank
  equal to $r\ind(\omega(p))$ whose endomorphism ring has Brauer class
  $\omega(p)$.  Let $\xi$ be the residual gerbe at $p$ (the fiber of
  $\ms M^{s, \mu}_{\ms X/k} \to M^{s, \mu}_{\ms X/k}$).  By definition,
  there is a $\xi$-twisted $\ms X$-twisted stable sheaf $\ms V$ on $\ms
  X\times\xi$ of rank $r$ with endomorphism ring $k$.  If $F$ is a
  $(-1)$-fold $\xi$-twisted vector space, then $F\tensor\ms V$ is a
  geometrically polystable $\ms X$-twisted sheaf with geometric
  components supported at $p$.  Furthermore, $\End(F\tensor\ms
  V)\cong\End(F)$ and this has Brauer class $\omega(p)$.  The rank of
  $F$ is a multiple of $\ind(\omega(p))$, and all multiples occur.
  Choosing one of minimal rank yields the desired result.  (In fact, we
  have shown that any $V$ has the form $F\tensor\ms V$ for some $F$.)
\end{proof}

\section{Index reduction via stable twisted sheaves}
\label{sec:index-reduction-via}

\subsection{Twisted Riemann-Roch}
\label{sec:twisted-riemann-roch}

We briefly summarize a twisted version of the Riemann-Roch theorem for
gerbes.  The reader is referred to section 2.2.7 of
\cite{twisted-moduli} for a proof.  Let $f:X\to S$ be a proper lci
morphism between quasi-projective varieties over a field.  Let $\ms
S\to S$ be a $\m_n$-gerbe and let $\mf f:\ms X\to\ms S$ be the
pullback of $\ms S$ to $X$.  Write $\td_f$ for the relative Todd class
of $f$ (the Todd class of the relative tangent complex).  Let $\ms F$
be a coherent $\ms X$-twisted sheaf.  For the purposes of this paper, when $n$ is invertible in the base field, 
define the rational Chow groups $A$ of $\ms S$ and $\ms X$ to be the
rational Chow groups of $S$ and $X$, respectively.  This is justified
by the results of Vistoli's thesis \cite{vistoli}; one can check that there is
a theory of Chern classes, etc., for this choice of Chow theory.  When $n$ is divisible by the characteristic of the base field, one must use Kresch's more general theory for Artin stacks \cite{kresch}.  Since the results we describe are all reduced to the case in which $n$ is invertible on the base, we need not concern ourselves with the subtleties of Kresch's Chow groups.

\begin{prop}\label{sec:twisted-riemann-roch-1}
  There is an equality $\chern(\mf f_{\ast}\ms F)=\mf
  f_{\ast}(\chern(\ms F)\cdot\td_f)$ in $A(\ms S)$.
\end{prop}

The case of primary interest to us will be when $S=\spec k$ and $X$ is a
curve.  In this case, there is an ad hoc proof of the following
corollary.

\begin{cor}\label{sec:twisted-riemann-roch-2}
  Let $f:X\to\spec k$ be a proper smooth curve and $\ms S\to\spec k$ a
  $\m_n$-gerbe with pullback $\ms X\to X$.  Write $\mf f:\ms X\to\ms
  S$ for the projection.  Given a coherent $\ms
  X$-twisted sheaf, the rank of the complex $\R\mf f_{\ast}\ms F$
  of $\ms S$-twisted vector spaces equals $\deg(\ms F)+\rk(\ms F)(1-g)$.
\end{cor}
\begin{proof}
  We give the proof in this case, as we will use it and it is simpler
  than the general case.  One need only note that the formation of
  $\R\mf f_{\ast}$ commutes with flat base change (Proposition 13.1.9
  of \cite{l-mb}), and that there is a finite extension $L$ of $k$ and
  a flat map $\spec L\to\ms S$.  Pulling back to $L$ reduces this to
  the classical Riemann-Roch formula.
\end{proof}
\begin{remark}
  The reasoning used in the proof of \ref{sec:twisted-riemann-roch-2}
  also proves \ref{sec:twisted-riemann-roch-1}, but one must pay more
  attention to the properties of the Chow theory.
\end{remark}

\subsection{A formula for index reduction}
\label{sec:form-index-reduct}

In this section $X$ is a smooth proper geometrically connected curve
over a perfect field $k$.  Write $D$ for the index of $X$ (the g.c.d. of
the degrees of all closed points) and $\delta$ for the index of
$\Pic^1_{X/k}$.  We clearly have that $\delta | D$ Let
$\widebar{\beta}\in\Br(k)$ be a Brauer class and $\beta\in\H^2(\spec
k,\m_n)$ a lift. 

In the following, given a scheme $Y$, the notation ``$p\in Y$'' will
mean that $p$ is a closed point of $Y$.  Given a Brauer class
$\alpha\in\Br(Y)$ and $p\in Y$, we will write $\alpha(p)$ for
$\alpha|_{\spec\kappa(p)}\in\Br(\kappa(p))$.  

\begin{defn}\label{sec:index-reduction-via-1}
  Given a scheme $Y$ and a Brauer class $\alpha\in\Br(Y)$, define the
  \emph{$\beta$-index reduction of $(Y,\alpha)$\/} by
  $$\iota_{\beta}(Y,\alpha)=\min_{p\in
  Y}[\kappa(p):k]\ind(\alpha(p)+\beta).$$
\end{defn}

Given $r$ and $d$, write $\iota_{\beta}(r,d)$ for
$\iota_{\beta}(M_{X/k}^s(r,d),\ms M^s_{X/k}(r,d))$.

\begin{prop}\label{sec:index-reduction-via-2}
  The index of $\beta_{k(X)}$ is $$\ind(\beta_{k(X)})=\min_{r|i,
    d\in[0,D)}r\iota_{\beta}(r,rd).$$ Furthermore, we have
  that $$\min_{d\in[0,D)}\iota_{\beta}(1,d)\text{ divides
  }\delta\ind(\beta_{k(X)}).$$
\end{prop}

The reader will note that the divisibility statement was originally
proven by Schofield and Van den Bergh \cite{schofield-vandenbergh}.
The present techniques give a new (but very closely related) proof of
the result.

Let us set notation for the proof. Choose a $\m_n$-gerbe $\ms S\to\spec
k$ representing $\beta$.  There is a $\m_n$-gerbe $\ms X\to X$ induced
by pullback of $\ms S$.  Write $i$ for the index of $\beta$ (over $k$).

\begin{proof}
  Given a locally free $\ms X$-twisted sheaf $V$ of minimal rank, we
  see that (1) $V$ is stable of some slope $\mu$, and (2) $V=S(V)$,
  the $\mu$-socle of $V$.  It follows that we may assume that $V$ is
  geometrically polystable.  By \ref{sec:moduli-stable-sheav-5} and
  \ref{sec:simple-bundl-pullb-4}, we see that the index of
  $\beta_{k(X)}$, which is the rank of $V$, must equal the
  $\beta$-index reduction of $(\M^s_{X/k},\ms M^s_{X/k})$ at $p$ for
  some closed point $p\in\M^s_{X/k}$.  It remains to show that the
  slope $\mu$ is integral, i.e., that $r | d$.  To see this, let
  $\pi:\ms X\to\ms S$ denote the natural morphism.  The Riemann-Roch
  formula \ref{sec:twisted-riemann-roch-2} shows that $\R\pi_{\ast}V$
  is a complex of $\ms S$-twisted vector spaces of rank $d+r(1-g)$.
  We know that $i$ must divide this quantity, and that $r$ must divide
  $i$ (as the index cannot increase upon pullback!).  By definition,
  $D$ is the minimal degree of a Cartier divisor on $X$.  It follows
  that twisting $V$ by an appropriate invertible sheaf allows us to
  assume that $d$ lies between $0$ and $D-1$, which proves the first
  statement.
  
  To prove the second, suppose there is a $k$-rational point $q$ on
  $\Pic^{\delta}_{X/k}$.  Thus, there is an $\alpha(q)$-twisted
  invertible sheaf $L$ of degree $\delta$.  Applying
  \ref{sec:twisted-riemann-roch-2} to the twisted sheaf $V\tensor
  L^{\tensor n}$ yields a complex of $\beta+n\alpha(q)$-twisted vector
  spaces of rank equal to $d+nr\delta+r(1-g)$.  The $\gcd$ of these
  ranks divides $r\delta$, which almost yields the desired statement.
  To see that the $\gcd$ can be replaced by $\min$ (as the formula in
  this case requires), we use an argument essentially due to Schofield
  and van den Bergh.  We may first replace $\beta$ by its $p$-primary
  component and assume that $\ind(\beta_{k(X)})$ is a power of
  $p$.  (We implicitly use the fact that the set of $k$-rational
  Picard obstructions is a group in reassembling the result from its
  primary parts.)  We may also replace $\delta$ by the largest
  $p$-power dividing it.  Writing $d=rd'$, the twisted Euler
  characteristic equals $r(d'+n\delta+1-g)$.  Choosing $n$
  appropriately, we see that $\ind(\beta+n\alpha(q))=r\delta'$ with
  $0\leq\delta'<\delta$, so that $\ind(\beta+n\alpha(q))<r\delta$.  It
  follows that the index of the $p$-primary part must divide
  $r\delta$, as required.  Furthermore, the $p$-primary part of
  $\beta+n\alpha(q)$ is $\beta+n\alpha(q)^{(p)}$, and it is easy to
  see that if $\alpha$ is a Picard obstruction then all of its primary
  parts are also Picard obstructions.  The result follows.
\end{proof}

This result has a particularly nice interpretation for curves of genus
$1$ (and $0$, although this case is already well known).

\begin{cor}\label{sec:index-reduction-via-3}
  If $g(X)\leq 1$ then $\ind(\beta_{k(X)}) = \min\{[E:k] \mid \beta_{X_E}
  \text{ is split}\}$.
\end{cor}

\begin{remark}
  We may interpret this in a few different ways. For one, it says that
  if $D$ is a central division algebra over $k(X)$ representing the
  Brauer class $\beta_{k(X)}$, then $D$ has a maximal subfield of the
  form $E \otimes_k k(X)$ for $E/k$ finite. From another perspective, it
  says that we may reduce the computation of index reduction to an
  understanding of splitting: $\beta$ has index $m$ if and only if there
  is a finite extension $E/k$ of degree $m$ such that the class
  $\beta_E$ is split by the curve $X_E$.
\end{remark}

We are able to derive from this a similar result in the case of an
imperfect field:
\begin{cor}\label{index-reduction-nonperfect}
  Let $k$ be an imperfect field of characteristic $p$ and suppose $p$ does not divide 
  $\ind(\beta)$. 
  If $g(X)\leq 1$ then $\ind(\beta_{k(X)}) = \gcd\{[E:k] \mid \beta_{X_E}
  \text{ is split}\}$.
\end{cor}
\begin{proof}
  Let $F/k$ be the perfect closure of $k$. In other words, $F$ is a
  perfect field lying in a fixed algebraic closure $\overline{k}$ of $k$
  which is closed under adjoining $p$th roots of elements, and is a
  compositum of $p$-power extensions. By the assumption on the
  characteristic, we may find a finite extension $E/F$ of degree $i =
  \ind(\beta_{X_F})$ such that $\beta_{X_E}$ is split. Since the
  condition of $\beta_{X_E}$ being split involves a finite set of
  equations with a finite number of elements of $E$, we may find a
  finitely generated $k$-subfield $k\subset E' \subset E$ such that $\beta_{X_{E'}}$
  is split. Since $E/k$ is algebraic, we have $[E' : k]$ is finite and
  must divide $i p^l$ for some $l \geq 0$. On the other hand, since $p$ does not divide $n$, $\beta$ is split by a finite extension $L/k$ of degree
  prime to $p$. Therefore $\gcd\{[E:k], [L:k]\} | i$, forcing the
  desired conclusion.
\end{proof}

Let us denote by $\alpha$ the Brauer class of the gerbe $\sPic_X \to
\Pic_X$. To warm up to the proof of Corollary
\ref{sec:index-reduction-via-3}, we give the following lemma (which is
actually a special case):
\begin{lem} \label{pic_splits} The class $\beta_{k(X)}$ is trivial if
  and only if $\beta = \alpha(p)$ for some point $p \in \Pic_X(k)$.
\end{lem}
\begin{proof}
  We note that $\beta_{k(X)}$ is trivial if and only if there is a
  $\beta$-twisted invertible sheaf on $X$, or in other words, if there
  is an object in $\sPic_{\ms X}(k)$. By \ref{C:tw-pic-obs-ch}, we have an
  identification $\Pic_{\ms X} = \Pic_X$ and using this, the Brauer
  class of the gerbe $\sPic_{\ms X} \to \Pic_{\ms X} \cong
  \Pic_X$ is given by $\alpha + \beta$. Therefore, noting that having an
  object of $\sPic_{\ms X}(k)$ is equivalent to having a point in $p \in
  \Pic_X(k)$ with trivial obstruction $\alpha(p) + \beta$, this says
  $\beta_{k(X)}$ is split if and only if there is a $p \in
  \Pic_X(k)$ with $\alpha(p) = -\beta$.

  To finish, we note that since the kernel $\Br(k) \to \Br(k(X))$ is a
  subgroup, $\beta_{k(X)}$ is split if and only if $-\beta_{k(X)}$ is
  split. The above argument then shows this to be equivalent to the
  existence of a $p \in \Pic_X(k)$ with $\alpha(p) = \beta$.
\end{proof}

\begin{proof}[Proof of Corollary \ref{sec:index-reduction-via-3}]
  Using the fact that any stable vector bundle on a genus $1$ curve with
  $r | d$ is invertible (\cite{atiyah}), and any stable vector bundle on
  a genus $0$ curve is invertible, it follows from
  \ref{sec:index-reduction-via-2} that we have
  \begin{align*}
    \ind(\beta_{k(X)}) &= \min_{d \in [0, \delta)}
    \iota_{\beta}(1,d) \\
    &= \min_d \big\{\min\{[k(p):k] \ind(\beta + \alpha(p)) \mid  p \in
    \Pic^d_X\}\big\} \\
    &= \min\{[k(p):k] \ind(\beta + \alpha(p))\mid p \in \Pic_X\} \\
    &= \min\{[L:k] \ind(\beta + \alpha(p))\mid p \in \Pic_X(L)\} \\
    &= \min\{[L:k] \ind(\beta_L + \gamma)\mid \gamma \in \Br(X_L/L)\}
    \text{ \ \ (by Lemma \ref{pic_splits})} \\
    &= \min\{[L:k] \ind(\beta_L - \gamma)\mid \gamma \in \Br(X_L/L)\} \\
    &= \min\{[L:k] [E:L]\mid \beta_E \sim \gamma_E, \gamma \in Br(X_L/L)\}
    \\
    &= \min\{[E:k]\mid \beta_{X_E} \text{ is split}\}
\end{align*}
\end{proof}

\subsection{Index reduction for local fields}

We can use Corollary \ref{sec:index-reduction-via-3} to give precise
information about index reduction in the case of local fields, using the
work of Roquette \cite{Roq} which computes
the relative Brauer group of a curve in this case.

Let $k$ be a local field, and let $C/k$ be a curve (which will soon be
assumed to have genus $1$). We recall the following result:

\begin{thm}[\cite{Roq}, Theorem 1]
  Let $A$ be a central simple $k$ algebra. Then $A \otimes_k k(C)$ is
  split if and only if $\ind(A) | \ind(C)$.
\end{thm}

By standard facts from local class field theory, we know that for $E/k$
finite, $$\ind(A_E) = \ind(A)/\gcd\{\ind(A), [E:k]\}.$$  Now suppose that
$C$ has genus $1$, and write $i = \ind(A)$. We may rewrite Corollary
\ref{sec:index-reduction-via-3} as saying:
\begin{equation*}
  \ind(A_{k(C)}) = \gcd\left\{[E:k]\ \left|\ \frac{i}{\gcd\{i, [E:k]\}}
  \big| \ind(C_E) \right.\right\}.
\end{equation*}

The computation of the index of $A_{k(C)}$ may therefore be expressed
entirely in terms of arithmetic information about the curve $C$, in
particular, how to compute its index over different finite extensions.
To give an example of this, let us consider the case where the index of
$C$ is $p$, a prime number, and suppose $\ind(A) = mp^n$, where $p \not|
m$, $n > 0$. For such a curve $C$, we define its capacity $\cpc(C)$ as:
$$\cpc(C) = \max\left\{r \left| 
\begin{matrix}
  \exists L/k \text{ a finite field extension with} \\ [L:k] = m'p^{r},
  p \not | m', \text{ and } C(L) = \emptyset 
\end{matrix}
\right.\right\}.$$
We then have:
$$\ind(A_{k(C)}) = \left\{
\begin{matrix}
  mp^n && \text{if $\cpc(C) < n-1$} \\ mp^{n-1} && \text{if $\cpc(C)
  \geq n-1$}
\end{matrix}
\right.$$

\subsection{Higher-dimensional varieties}
\label{sec:high-dimens-vari}

We indicate in this section how to extend our results to higher
dimensional varieties over $k$.  For the most part, the results are
straightforward generalizations of the techniques above.  One must do
slightly more numerical work with the Riemann-Roch formula (as in
\cite{schofield-vandenbergh}).  There is also a slight complcation coming
from the
difference between the category of torsion free sheaves and its 
quotient by the category of sheaves supported in codimension at least
$2$ (which is necessary to mimic the argument reducing to the
geometrically polystable case).  
We are content to simply state the results and leave the
mostly straightforward details to the reader.

Let $X$ be a smooth geometrically connected projective variety over $k$
of dimension $t$ with fixed ample invertible sheaf $\ms O(1)$.  Any
section $\sigma$ of $\Pic_{X/k}$ has a well-defined degree given by the
top self-intersection of a divisor on $X\tensor\widebar k$ representing
$\sigma$.  Since the degree is an intersection-theoretic invariant, it
is clearly constant on connected components of $\Pic_{X/k}$.  We will
write $\sPic_{X/k}^d$ for the stack of invertible sheaves of degree $d$;
the usual results show that $\sPic^d_{X/k}$ is a $\G_m$-gerbe over its
sheafification $\Pic_{X/k}^d$.

Let $\delta$ be the $\gcd$ of the degrees of all $k$-rational sections
of $\Pic_{X/k}$ and let $D$ be the $\gcd$ of the degrees of all
$k$-rational sections of $\sPic_{X/k}$ (i.e., those sections of
$\Pic_{X/k}$ which arise from actual invertible sheaves on $X$).

\begin{defn}
   Given a torsion free coherent sheaf $\ms E$ of rank $r$ on $X$, the
   \emph{generalized slope\/} of $\ms E$ is $\chi(\ms E)/r$.
\end{defn}

Recall that a torsion free sheaf $\ms E$ of positive rank $r$ on $X$
is \emph{slope semistable (resp.\ stable) with respect to $\ms
  O(1)$\/} if for all proper subsheaves $\ms F\subset\ms E$ with
strictly smaller rank, we have that $c_1(\ms F)/\rk(\ms F)-c_1(\ms
E)/r$ has non-positive (resp.\ strictly negative) intersection with
$\dim(X)-1$ copies of a divisor in $|\ms O(1)|$.  It is a standard
result that the stack of slope semistable sheaves is an Artin stack
with the substack of slope stable sheaves an open substack which is a
$\G_m$-gerbe over an algebraic space.  We will denote the stack of slope
stable torsion free sheaves of rank $r$ and generalized slope $\mu$ by 
by $\ms M^s_{X/k}(r,\mu)$ and its sheafification by $M^s_{X/k}(r,\mu)$.

Let $\beta\in\Br(k)$ have index $i$.  Let $\ms S\to\spec k$ be a
$\m_n$-gerbe representing $\beta$ and $\ms X=X\times_k\ms S$.  We will
write $\mf f:\ms X\to\ms S$ for the natural map.  Given an $\ms
X$-twisted sheaf $\ms F$, we will write $\chi(\ms F)$ for the rank of
$\R\mf f_{\ast}\ms F$ (as a complex of $\ms S$-twisted sheaves).

\begin{prop}\label{P:ind-red-higher-dim}
  The index of $\beta_{k(X)}$ is
$$\ind(\beta_{k(X)})=\min_{r | i,\mu\in\Z}r\iota_{\beta}(M^s_{X/k}(r,\mu),\ms M^s_{X/k}(r,\mu)).$$
Furthermore, we have that
$$\min_{d\in[0,D)}\iota_{\beta}(\Pic_{X/k}^d,\sPic_{X/k}^d)\text{
  divides }\delta\ind(\beta_{k(X)}).$$
\end{prop}
As above, the second statement is (for higher dimensional varieties) a
refinement of the main result of Schofield and Van den Bergh: they
consider only very ample sections of $\Pic$ when computing $\delta$,
while an argument in the derived category as above shows that in fact
one can strengthen the result (i.e., lower $\delta$) by considering
arbitrary sections of $\Pic$.

The proof of \ref{P:ind-red-higher-dim} follows precisely the outline of
the proof of \ref{sec:index-reduction-via-2} above.  We sketch the proof
of the second statement, following p.\ 732 of
\cite{schofield-vandenbergh} essentially verbatim as an aid to the
reader.  We may assume that $\beta$ is $p$-primary for some prime $p$,
so that we may assume that $\delta$ is realized by a $k$-point $q$ of of
$\Pic_{X/k}$, corresponding to some $\alpha(q)$-twisted invertible sheaf
$\ms L$.  Let $V$ be a torsion free $\ms X$-twisted sheaf of rank
$\ind(\beta_{k(X)})$.  It follows from \ref{sec:twisted-riemann-roch-1}
that
$$\chi(V\tensor\ms L^{\tensor
m})=\frac{\delta\ind(\beta_{k(X)})}{t!}m^t+\text{ lower order terms.}$$
Standard manipulations of numerical polynomials show that for a fixed
$m$, we have
$$\delta\ind(\beta_{k(X)})=\sum_{j=0}^t(-1)^j\binom{t}{j}\chi(V\tensor\ms
L^{m+t-j}).$$
Since $\delta=p^a$ and $\ind(\beta_{k(X)})=p^b$ for some $a$ and $b$, it
follows that there is some $N$ such that $\chi(V\tensor\ms L^{N})$ is
divisible by at most $p^{a+b}$.  This is easily seen to imply the
$p$-primary part of the required statement.  The general case follows by
reassembling the primary parts.

\section{Index reduction on torsors under abelian varieties and homogeneous twisted bundles}
\label{sec:index-reduct-abel}

In this section we sketch an approach to index reduction on torsors
under abelian varieties using twisted Fourier-Mukai transforms. As we
show, there is a connection between the existence of homogeneous
minimal bundles and an index reduction formula involving only the
universal Picard obstruction.  Such homogeneity results are perhaps
analogous to the explicit vector bundles used in the index
reduction formulas arising in \cite{Panin,Merk-Panin-Wadsworth} for
certain homogeneous spaces under linear algebraic groups using the
$K$-theory of such spaces.  The main result of this section may be
interpreted as evidence that homogeneous index reduction is
exceedingly unlikely for torsors under abelian varieties.  Since the
$K$-theory of such a torsor is not generated by equivariant sheaves
(in contrast to the rational case), this should not be surprising.

\subsection{Twisted Fourier-Mukai transforms}
\label{sec:twist-four-mukai}

In this section we prove a twisted form of Mukai's theorem on derived
equivalences of abelian varieties, relating twisted sheaves on an
$A$-torsor to twisted sheaves on $A$ in a certain Brauer class.  For the
sake of simplicity, we restrict our attention to bounded derived
categories of twisted sheaves.

Let $k$ be a field and $A$ an abelian variety over $k$ and $T$ an
$A^\vee$-torsor.  Suppose $\beta\in\Br(k)$ is a Brauer class and let
$\ms T\to T$ be a $\G_m$-gerbe representing $\beta_T$.

\begin{defn}
  The \emph{twisted Picard stack\/} parametrizing invertible $\ms T$-twisted
  sheaves will be denoted $\sPic_{\ms T/k}$.
\end{defn}

Applying \ref{C:tw-pic-obs-ch}, we see that $\sPic_{\ms T/k}$ is a
$\G_m$-gerbe over $\Pic_{T/k}$ with Brauer class
$[\sPic_{T/k}]+\beta$.  In particular, there is a distinguished
connected component $\sPic^0_{\ms T/k}$ corresponding to the component
of $\sPic_{T/k}$ containing the point $\ms O_T$.  It follows that
$\sPic^0_{\ms T/k}$ is a $\G_m$-gerbe over $A$, such that the fiber
over the identity section, viewed as a $\G_m$-gerbe, has Brauer class
$\beta$.  We will write $\ms A\to A$ for this gerbe in what follows
(the class $\beta$ being understood throughout).

There is a universal invertible sheaf $\ms L$ on $\ms A\times\ms T$;
the geometric fibers of $\ms L$ over $\ms T$ are $\ms A$-twisted,
while the geometric fibers over $\ms A$ are $\ms T$-twisted.  
We can thus define a Fourier-Mukai transform $\Phi:\D_{-1}^b(\ms
A)\to\D_1^b(\ms T)$ with kernel $\ms L$ from the derived category of $(-1)$-fold
$\ms A$-twisted coherent sheaves to the derived category of coherent
$\ms T$-twisted sheaves.

\begin{prop}\label{P:FM-equiv}
  The functor $\Phi$ is an equivalence of categories.
\end{prop}
\begin{proof}
  By an argument formally identical to Lemma 2.12 of \cite{orlov}, it
  suffices to prove the statement when $k$ is algebraically closed.
  In this case, $\ms T\to T$ and $\ms A\to A$ are trivial gerbes.
  Given a trivialization of $\ms T\to T$ and a $k$-point of $T$, one
  naturally gets an identification $T\simto A^{\vee}$ and a
  trivialization of $\ms A\to A$ (e.g., by the standard method of
  ``rigidifications'' of invertible sheaves \cite{blr}).  Moreover,
  trivializations of $\ms T\to T$ and $\ms A\to A$ serve to identify
  the derived categories of $n$-fold $\ms T$-twisted (resp.\ $\ms
  A$-twisted) sheaves with coherent sheaves on $T$ (resp.\ $A$) for
  any $n$.  It is easy to check that composing $\Phi$ with these
  equivalences yields the standard Fourier-Mukai transform
  $\D^b(A)\to\D^b(A^{\vee})$, which is an equivalence by Theorem 2.2
  of \cite{mukai}.
\end{proof}

\subsection{Moduli of homogeneous bundles}
\label{sec:moduli-homog-bundl}

Using the twisted Fourier-Mukai transform, we can describe certain
moduli spaces of locally free twisted sheaves on $T$ in terms of moduli spaces
of finite length twisted sheaves on $\ms A$.

Let $\beta\in\Br(k)$, and let $\ms A\to A$ and $\ms T\to T$ be as
above.  Note that in constructing $\ms T\to T$ we may explicitly use
the $1$-fibered product of $T$ with a $\G_m$-gerbe over $\spec k$
representing $\beta$.  It follows that $A$ acts on $\ms T$
functorially (and not merely pseudo-functorially).  This allows us to
think about the pullback of a $\ms T$-twisted sheaf via a translation
by a point of $A$ in a very concrete manner; we will implicitly do
this in what follows.

We will write $p$ for the first projection $A\times T\to A$, $q$ for
the second projection $A\times T\to T$, and $\mu$ for the action
$A\times T\to T$.

\begin{defn}
  A $\ms T$-twisted sheaf $\ms F$ is 
  \begin{enumerate}
  \item \emph{homogeneous\/} if $\tau_x^{\ast}\ms
    F_{\kappa(x)}\cong\ms F_{\kappa(x)}$ for all geometric points
    $x:\spec\kappa\to A$;
  \item \emph{uniformly homogeneous\/} if the scheme
    $\isom_A(\mu^{\ast}\ms F,q^{\ast}\ms F)\to A$ has sections
    \'etale-locally on $A$;
  \item \emph{semi-homogeneous\/} if for all geometric points
    $x:\spec\kappa\to A$, there exists an invertible sheaf $\ms
    L_x\in\Pic(T_{\kappa})$ such that $\tau_x^\ast\ms
    F_{\kappa}\cong\ms L_x\tensor\ms F_\kappa$ (where, by abuse of
    notation, we also let $\ms L_x$ denote the pullback of $\ms L_x$
    to $\ms T_\kappa$).
  \end{enumerate}
\end{defn}

Note that uniformly homogeneous sheaves are automatically homogeneous,
and that homogeneous sheaves are locally free.  The classical
literature deals with abelian varieties over algebraically closed
fields, and the usual definition of homogeneity works only with
geometric points with coefficients in the base field.

\begin{lem}
  A $\ms T$-twisted sheaf $\ms F$ is homogeneous if and only if
  $\tau_x^{\ast}\ms F_{\kappa(x)}\cong\ms F_{\kappa(x)}$ for
  arbitrary field-valued (not necessarily geometric) points
  $x:\spec\kappa\to A$.
\end{lem}
\begin{proof} Suppose $\ms F$ is homogeneous.  Without loss of
  generality, we may assume $\kappa(x)=k$.  By the compatibility of
  cohomology with flat base change, the sheaf of isomorphisms
  $\Isom(\tau_x^{\ast}\ms F,\ms F)$ is represented by an open subscheme
  of the affine space $\hom(\tau_x^{\ast}\ms F,\ms F)$; this subscheme
  is non-empty by the homogeneity assumption.  Thus, if $k$ is infinite,
  it immediately follows that there is a rational point.  If $k$ is
  finite, then we need only note that the $\Isom$ scheme is a torsor
  under the automorphism sheaf $\saut(\ms F)$, which is a smooth
  geometrically connected group scheme over $k$ (see
  \ref{sec:moduli-homog-bundl-1} below).  By Lang's theorem, any torsor
  has a rational point.
\end{proof}

\begin{cor}
  Given a field extension $K/k$, a $\ms T$-twisted sheaf $\ms F$ is
  homogeneous if and only if $\ms F\tensor K$ is homogeneous.
\end{cor}

Given two sheaves $F$ and $G$ on $T\times A$, define a functor
$\uhom_A(F,G)$ on the category of $A$ schemes by sending $\gamma:S\to
A$ to $\hom_S(F_S,G_S)$.

\begin{lem}\label{sec:moduli-homog-bundl-1}
  The scheme $\isom_A(\mu^{\ast}\ms F,q^{\ast}\ms F)\to A$ The functor
  $\uhom_A(\mu^{\ast}\ms F,q^{\ast}\ms F)$ is represented by a closed
  cone (with linear geometric fibers) in a geometric vector bundle and
  contains $\isom_A(\mu^{\ast}\ms F,q^{\ast}\ms F)$ as an open
  subscheme.
\end{lem}
\begin{proof}
  Let $\ms W$ be a locally free $\ms T$-twisted sheaf and $\ms O(1)$ an
  ample invertible sheaf on $T$.  Choose $N$ and $m$ such that there is
  a surjection $\pi:\ms W(N)^m\to \mu^{\ast}\ms F$ and such that
  $\R^ip_{\ast}\shom(\ms W,q^{\ast}\ms F)(N)=0$ for $i>0$.  Writing $\ms
  K$ for the kernel of $\pi$, we may similarly choose $N'\geq N$ and
  $m'$ such that there is a surjection $\ms W(N')^{m'}\to\ms K$.  There
  results a complex
  $$\ms W(N')^{m'}\to\ms W(N)^{m}\to \mu^{\ast}\ms F\to 0.$$
  The cohomological assumptions show that there is a geometric vector
  bundle $\mathbf V$ (resp.\ $\mathbf V'$) whose sheaf of sections is
  $$p_{\ast}(\shom(\ms W,\mu^{\ast}\ms F)(N)^m)$$ (resp.\
  $p_{\ast}(\shom(\ms W,\mu^{\ast}\ms F)(N')^{m'})$).  Moreover, the
  presentation of $\mu^{\ast}\ms F$ yields an exact sequence of abelian sheaves 
  $$0\to\uhom_A(\mu^{\ast}\ms F,q^{\ast}\ms F)\to\mathbf V\to\mathbf V',$$
  where the latter map is linear.  It follows that $\hom_A(\mu^{\ast}\ms
  F,q^{\ast}\ms F)$ is represented by a closed cone in $\mathbf V$ with
  linear fibers.  That isomorphisms form an open subscheme is immediate.
\end{proof}
  
\begin{prop}\label{sec:moduli-homog-bundl-2}
  A $\ms T$-twisted sheaf $\ms F$ is homogeneous if and only if it is
  uniformly homogeneous.  
\end{prop}
\begin{proof}
  Suppose $\ms F$ is homogeneous, so that the open subscheme
  $$\isom_A(\mu^{\ast}\ms F,q^{\ast}\ms F)\subset\uhom_A(\mu^{\ast}\ms
  F,q^{\ast}\ms F)$$ meets every fiber.  In this case, since
  isomorphisms are dense in every fiber (the fibers being geometrically
  integral) and each fiber is a torsor under $\aut(\ms F)$, we see that
  the fibers of $\uhom_A(\mu^{\ast}\ms F,q^{\ast}\ms F)$ are all of the
  same dimension.  Thus, the linear map $\mathbf V\to\mathbf V'$ from
  the proof of \ref{sec:moduli-homog-bundl-1} has constant rank in every
  geometric fiber.  Since $A$ is reduced, we conclude that the kernel is
  a subbundle of $\mathbf V$, i.e., that $\uhom_A(\mu^{\ast}\ms
  F,q^{\ast}\ms F)$ is represented by a geometric vector bundle, hence
  is smooth over $A$.  As an open subscheme hitting every fiber,
  $\isom_A(\mu^{\ast}\ms F,q^{\ast}\ms F)$ must have sections everywhere
  \'etale-locally over $A$.
  
  The other direction is immediate.
\end{proof}
  
  
\begin{defn}
  A \emph{flat family of homogeneous $\ms T$-twisted coherent sheaves
    parametrized by $S$\/} is an $S$-flat quasi-coherent $\ms
  T$-twisted sheaf of finite presentation $\ms F$ on $\ms T\times S$
  such that for all geometric points $s\to S$, the fiber $\ms F_s$ is
  homogeneous.
\end{defn}

It is clear that the collection of flat families of homogeneous
coherent $\ms T$-twisted sheaves forms a stack.  We will prove in a
moment that it is in fact a quasi-proper Artin stack, although this is
not immediately obvious.  It is quite easy to see that it is a
constructible substack of the stack of coherent $\ms T$-twisted
sheaves (left to the reader), but we are not even sure if it is
locally closed.  (It is unlikely to be closed, as any homogeneous
sheaf is locally free.)  

We beg the reader's indulgence in allowing us the following notational
convenience.
\begin{con}
  Given a scheme $S$ and a stack $\ms S$ of categories on the fppf site of
  $S$, we will say that ``$\ms S$ is an Artin stack'' (etc.) if the
  underlying stack of groupoids is an Artin stack.  Thus, when we talk
  about stacks of sheaves, we include all morphisms in the category
  structure of the fiber categories, throwing away all but the
  isomorphisms only when we are considering the algebraicity
  properties of the stack in question.
\end{con}

\begin{defn}
  The stack of homogeneous $\ms T$-twisted sheaves of rank $r$ will be
  denoted $\HH^r_{\ms T/k}$.
\end{defn}

There is another stack which will be of interest in this section.  

\begin{defn}
  Given an Artinian $\ms A\tensor\widebar k$-twisted coherent sheaf
  $\ms F$, the \emph{length\/} of $\ms F$ is defined as follows: $\ms
  F$ admits a filtration $\ms F=\ms F^n\supset\ms
  F^{n-1}\supset\cdots\supset\ms F^0=0$ by sheaves $\ms F^i/\ms
  F^{i-1}$ with support equal to the reduced structures on residual
  gerbes of $\ms A$.  The length of $\ms F^i/\ms F^{i-1}$ is defined
  to be its rank as a sheaf on the residual gerbe and the length of
  $\ms F$ is then defined to be $\ell(\ms F)=\sum_i\ell(\ms F^i/\ms
  F^{i-1})$.
\end{defn}

While this definition may seem contrived, it fits naturally into a
general theory of Chern classes and Hilbert polynomials for twisted
sheaves.  In section 2.2.7 of \cite{twisted-moduli}, the reader will
find a proof that the length is constant in a flat family.  Whenever
we speak of the length of an $\ms A$-twisted coherent sheaf $\ms G$,
we mean the length of $\ms G\tensor\widebar k$ (i.e., the length of
the geometric fiber).

\begin{defn}
  A \emph{flat family of finite length coherent $\ms A$-twisted
    sheaves parametrized by $S$\/} is an $S$-flat quasi-coherent $\ms
  A$-twisted sheaf $\ms F$ of finite presentation over $\ms A\times S$
  such that for every geometric point $s\to S$, the fiber $\ms F_s$ is
  a coherent $\ms A_s$-twisted sheaf of finite length.
\end{defn}

It is easy to see that the collection of finite length coherent $\ms
A$-twisted sheaves of finite length forms an Artin stack.  For
details, the reader can consult \cite{coherent-algebras} or section
2.3 of \cite{twisted-moduli}.

\begin{defn}
  The stack of finite length coherent $\ms A$-twisted sheaves of
  length $\ell$ will be denoted $\FF^\ell_{\ms A/k}$.
\end{defn}

We recall a well-known result of Mukai (which carries over 
\emph{mutatis mutandi\/} to the twisted case).

\begin{prop}[Mukai]
  For every algebraically closed extension field $K/k$, the
  Fourier-Mukai transform defines an equivalence of fiber categories
  $$\FF^r_{\ms A/k}(K)\simto\HH^r_{\ms T/k}(K).$$
\end{prop}
For the proof, the reader is referred to \S{3} of \cite{mukai}.  For
geometric purposes, it is useful to have the following totally
unsurprising extension of the classical result.

\begin{prop}
  Given any $k$-scheme $S$, the twisted Fourier-Mukai transform
  establishes an equivalence of categories $$\FF^r_{\ms
    A/k}(S)\simto\HH^r_{\ms T/k}(S).$$
\end{prop}
\begin{proof}
  By the obvious compatibility of the Fourier-Mukai transform with
  derived base change, it is clear that $\Phi$ yields functors between
  the categories described in the proposition.  Furthermore, it is
  clear that it suffices to prove the result assuming that $T$ has a
  section: any map of stacks which is an isomorphism locally on the
  base must be an isomorphism.  (Note: here we do not mean simply
  stacks in groupoids, but arbitrary stacks of categories.)  Arguing
  as in \ref{P:FM-equiv}, the statement reduces to the obvious
  (twisted) relativization of Theorem 2.2 of \cite{mukai}.
\end{proof}

\begin{cor}\label{C:isom-of-thingos}
  The Fourier-Mukai transform induces an isomorphism of stacks
  $\FF^r_{\ms A/k}\simto\HH^r_{\ms T/k}$.  In particular, $\HH^r_{\ms
    T/k}$ is a quasi-proper Artin stack of finite presentation over
  the base.
\end{cor}

\begin{remark}
  This gives a nice example of a proof of (non-obvious) algebraicity
  of a very concrete stack of coherent sheaves using methods which
  pass through the (not as concrete) fibered category of derived categories.
\end{remark}

\begin{ques}
  Is the stack of homogeneous coherent sheaves a locally closed
  substack of the stack of coherent sheaves?
\end{ques}

\subsection{A criterion for homogeneous index reduction}
\label{sec:appl-index-reduct}

Let $k$ be a field, $A$ an abelian variety over $k$, $T$ an
$A$-torsor, and $\beta\in\Br(k)$ a Brauer class.  Choose a gerbe $\ms
S\to\spec k$ representing $\beta$ and let $f:\ms T\to\ms S$ denote the
pullback to $T$.

\begin{defn}\label{D:homog-ind-red}
  The class $\beta$ has \emph{(semi-)homogeneous index reduction\/} if
  there is a complex of (semi-)homogeneous $\beta$-twisted sheaves on
  $T$ of rank equal to $\ind_{k(T)}(\beta_{k(T)})$.
\end{defn}

\begin{defn}\label{D:deg-0-obs}
  The class $\beta$ has \emph{index reduction by degree $0$
    obstructions\/} if
  $$\ind(\beta_{k(T)})=\iota_{\beta}(\Pic^0_{T/k},\sPic^0_{T/k}).$$
\end{defn}

The main result of this section is the following.

\begin{prop}
 The class $\beta$ has homogeneous index reduction if and only if it
 has index reduction by degree $0$ obstructions.
\end{prop}
\begin{proof}
  This follows immediately from \ref{C:isom-of-thingos}, which shows
  that the minimal rank of a homogeneous $\ms T$-twisted sheaf equals
  the minimal length of an Artinian $\ms A$-twisted coherent sheaf.
  But the latter is computed precisely by the formula given in
  \ref{D:deg-0-obs}.
\end{proof}

In the case of genus $1$ curves, example \ref{example:kra-mir} shows
that that homogeneous index reduction need not hold in general. However,
it follows immediately from \ref{sec:index-reduction-via-3} that any
genus $1$ curve admits semi-homogeneous index reduction.

\begin{example} \label{example:kra-mir}
Let $C/k$ be a genus $1$ curve with $\per C \neq \ind C$ (for the
existence of these, see Theorem 3 of \cite{pete} or \cite{cassels}). By
\cite{KraMir} Theorem 2.1.1, there is a Brauer class $\beta \in \Br(k)$
such that $\beta_{k(C)} = 0$, however $\beta$ is not equal to any
obstruction class from $\Pic^0_{E/k}$. In particular, $\beta$ does not
have index reduction by degree $0$ obstructions. 

It follows more generally that if $T/k$ is a torsor for an abelian
variety, a necessary condition that all Brauer classes
have index reduction by degree $0$ obstructions is that the divisorial
index group and the period group of $T$ coincide (see \cite{KraMir}, for
definitions of these groups and page 5 for the relevant conclusion). On
the other hand, it is not at all clear that the reverse implication
should hold and that for such torsors we should have homogeneous index
reduction.
\end{example}

\appendix

\section{A period-index result}
\vskip 1\baselineskip
\begin{center}
{\sc Bhargav Bhatt}
\end{center}
\vskip 1\baselineskip

Let $A \to S$ be an abelian scheme, and let $X$ be a torsor for it.
Our goal is to relate the derived category of coherent sheaves on $X$
to that on a suitable $\G_m$-gerbe on $A^t$, the dual abelian scheme.
As an application, we obtain a period-index result over global fields
for geometrically trivial Brauer classes on abelian varieties coming
from locally trivial torsors.

\subsection{Construction of the equivalence}

Let $\pi:A \to S$ be an abelian scheme, and let $\pi:A^t \to S$ be the
dual scheme. The Leray spectral sequence in the fppf topology for the latter with
$\G_m$-coeffecients gives us a low degree short exact sequence
$$1 \to \H^2(S,\G_m) \to F^1(\H^2(A^t,\G_m)) \to \H^1(S,\Pic_{A^t/S}) \to 1 $$
where $F^\bullet$ is the filtration defined by the spectral sequence. Note
that $F^1(\H^2(A^t,\G_m))$ can be identified as the subgroup of
$\H^2(A^t,\G_m)$ consisting of geometrically trivial classes (i.e:
classes killed by fppf localisation on $S$). The preceeding exact
sequence is canonically split by the zero section $e^t:S \to A^t$. In
particular, there's a natural (composite) map 
$$\alpha:\H^1(S,A) \cong \H^1(S,\Pic^0_{A^t/S}) \to \H^1(S,\Pic_{A^t/S})
\to \H^2(A^t,\G_m)$$
where $A \cong \Pic^0_{A^t/S} \to \Pic_{A^t/S}$ is the connected
component of the identity. Denote by $\alpha_X$ the Brauer class
associated to the class $[X]$ of a torsor $X$ under this map. Then there
are two ways interpret the association $X \mapsto \alpha_X$.

$\bullet$ \textit{A geometric interpretation}\\
Given an $A$-torsor $X$, we have a canonical associated $\G_m$-gerbe
$\mathbf{Pic}_{X/S} \to \Pic_{X/S}$. Restricting this $\G_m$-gerbe along
the connected component $\Pic^0_{X/S} \to \Pic_{X/S}$, we obtain a
$\G_m$-gerbe over $\Pic^0_{X/S}$. Since $X$ is an $A$-torsor, we have a
canonical isomorphism $A^t \cong \Pic^0_{X/S}$. Therefore, we obtain a
$\G_m$-gerbe $\mathbf{Pic}^0_{X/S} \to A^t$ whose class in
$\H^2(A^t,\G_m)$ will be denoted by $\beta_X$. Note that since a torsor
is geometrically trival, the functoriality of this association implies
that the class we obtain lands in $F^1(\H^2(A^t,\G_m))$.

$\bullet$ \textit{A homological interpretation}\\
Given an $A$-torsor $X$, the Leray spectral sequence for $X \to S$ (with
$\G_m$-coeffecients) gives us a low degree exact sequence
$$\Pic(X) \to \Pic_{X/S}(S) \to \H^2(S,\G_m)  $$
This sequence, of course, continues to exist if we base change on $S$.
Therefore, we can replace $S$ with $A^t \cong \Pic^0_{X/S}$, and $X$
with $X \times_S \Pic^0_{X/S}$ to get
$$\Pic(X \times_S \Pic^0_{X/S}) \to \Pic_{X/S}(\Pic^0_{X/S}) \to
\H^2(\Pic^0_{X/S},\G_m)$$
The natural inclusion $\Pic^0_{X/S} \to \Pic_{X/S}$ defines a canonical
element $\mathcal{P}$ (the ``Poincare bundle'') in
$\Pic_{X/S}(\Pic^0_{X/S})$. The obstruction to lifting $\mathcal{P}$ to
a line bundle on $X \times_S \Pic^0_{X/S}$ (i.e: the failure of
$\Pic^0_{X/S}$ to be a fine moduli space for translation invariant line
bundles on $X$) defines an element $\gamma_X$ in $\H^2(\Pic^0_{X/S},\G_m)
= \H^2(A^t,\G_m)$. Once again, by functoriality of the construction, it's
clear that the constructed element lies in $F^1(\H^2(A^t,\G_m))$.

A messy cocycle calculation shows that $\alpha_X = -\beta_X =
-\gamma_X$.  Given this, we construct the desired equivalence.
Specifically, we have:

\begin{btheorem}
\label{main thm}
Let $\pi:A \to S$ be an abelian scheme and let $\pi^t:A^t \to S$ be its
dual scheme. Given an $A$-torsor $f:X \to S$, let $\alpha_X \in
\H^2(A^t,\G_m)$ be the class defined above. Then the universal line
bundle $\mathcal{P}$ on $X \times_S \mathbf{Pic}^0_{X/S}$, when thought
of as the kernel of a Fourier-Mukai transform, defines an equivalence of
categories $\F:D(X) \to D(A^t,-\alpha_X)$
\end{btheorem}

\begin{proof}
This is proven in Theorem \ref{P:FM-equiv}. We'll simply remark that the
significance of the negative sign in $D(A^t,-\alpha_X)$ is that the
pushforward of $\mathcal{P}$ along the projection $X \times_S
\mathbf{Pic}^0_{X/S} \to \mathbf{Pic}^0_{X/S}$ is a $\gamma_X$-twisted
sheaf when the target is thought of as a $\G_m$-gerbe over $A^t$, and
that $\gamma_X = -\alpha_X$.
\end{proof}

\subsection{An application}

Our goal is to prove a period-index result for Brauer classes associated
to torsors for abelian varieties. Specifically, we will show:

\begin{btheorem}
\label{main app}
Let $A$ be a $g$-dimensional principally polarised abelian variety over
a field $k$. Given an $A$-torsor $X$, let $\alpha_X \in \H^2(A^t,\G_m)$
be the class defined earlier. Assume one of the following
\begin{enumerate}
\item $k$ has trivial Brauer group.
\item $k$ is a global field and $X$ is an element of $\Sha(A)$ (i.e: has
points locally).
\end{enumerate}
Then $\ind(\alpha_X) | \per(X)^g$. Furthermore, if $X$ has odd order
in in $\H^1(k,A)$, then $\ind(\alpha_X) | \per(\alpha_X)^g$.
\end{btheorem}

In the above situation, Theorem 4 of \cite{PeteClark} proves that 
$\ind(X) \leq g! \per(X)$. Even without any assumptions on $k$, an
argument due to Lenstra (explained in Proposition 12 of \cite{PeteClark})
shows that $\ind(\alpha_X) \mid \per(X)^{2g}$ provided that $\per(X)$ is 
invertible in $k$. Thus, one can view Theorem \ref{main app} as an
improvement of both these results when interpreted in terms of the Brauer
group. On the other hand, in the case that $k$ is the function field of 
a curve over $\C$, it has been conjectured (see page12 of \cite{ctl}) that 
one should have $\ind(\alpha_X) | \per(\alpha_X)^{g}$ for any Brauer class
on a variety of dimension $g$. By Tsen's theorem, we may apply Theorem 
\ref{main app} and verify that the conjecture holds in the case of abelian
varieties.

The idea of the proof is to use Theorem \ref{main thm} to construct an
$\alpha_X$-twisted vector bundle of rank $\per(X)^g$ on $A^t$. In order
to do so, however, we need a good supply of sheaves on $X$ which is exactly
what the following Lemma (proven differently in Proposition 21 of 
\cite{PeteClark}) accomplishes.

\begin{blemma}
\label{descent lemma}
Assume we're in the situation of Theorem \ref{main app}. If $L \in
\Pic(A)$ is a line bundle defining a principal polarisation, then the
class of $L^{\per(X)}$ in $\NS_{A} \cong \NS_{X}$ is effective
(representable by a line bundle on $X$).
\end{blemma}
\begin{proof}
We'll think of the torsor $X$ as an extension
$$1 \to A \to T \to \Z \to 1 $$
where the fibre of $T \to \Z$ over $1 \in \Z$ is the torsor $X$. Given a
line bundle $L \in \Pic(A)$ which defines a principal polarisation
$\phi_L:A \to A^t$, we can push this exact sequence out along $\phi_L$
to obtain an extension
\begin{equation}
\label{ext 1}
1 \to A^t \to T' \to \Z \to 1 
\end{equation}
which corresponds to the image of $X$ under $\H^1(\phi_L):\H^1(k,A) \to
\H^1(k,A^t)$. As a class in $\Ext^1_{k}(\Z,A^t)$, this has order
$\per(X)$. On the other hand, since $X$ is an $A$-torsor, we have
canonical isomorphisms $A^t \cong \Pic^0_{X}$ and $\NS_{A} \cong
\NS_{X}$. Thus, we obtain a canonical extension 
$$1 \to A^t \to \Pic_{X} \to \NS_{A} \to 1 $$
The class of the line bundle $L$ in $\NS_{A}$ defines a morphism $\Z \to
\NS_{A}$. Pulling the preceeding exact sequence back along this
morphism, we obtain an extension
\begin{equation}
\label{ext 2}
1 \to A^t \to T'' \to \Z \to 1 
\end{equation}

A messy calculation with cocyles reveals that the class of this
extension in $\Ext^1_k(\Z,A^t)$ is the inverse of the class associated
to extension \eqref{ext 1}. The upshot of this is that extension
\eqref{ext 2} also has order $\per(X)$ in the group of extension
classes. This means that there is an element $M \in \Pic_{X}(k)$ whose
class $\NS_{X} \cong \NS_{A}$ is the same as that of $L^{\per(X)}$. What
remains to be shown now is that $M \in \Pic_{X}(k)$ comes from an actual
line bundle on $X$. The obstruction to $M$ being represented by a line
bundle lives in $\Br(k)$. In the case that $\Br(k) = 0$ we are obviously
done. In the case of global fields, recall that an element of $\Br(k)$
is trivial if it is locally trivial. By functoriality of the obstruction
and the fact that $X$ is locally split, we see that there is no
obstruction locally. Thus, there is no global obstruction to
representing $M$ by an actual line bundle in this case as well.
\end{proof}

We're now in a position to give a proof of Theorem \ref{main app}

\begin{proof}[Proof of Theorem \ref{main app}]
Note that the kernel of $$\alpha:\H^1(k,A) \to
\H^1(k,\Pic_{A^t}) \to \H^2(A^t,\G_m)$$ is $2$-torsion because the kernel
of the first map is $2$-torsion\footnote{There is a short exact sequence
$\H^0(k,\NS_{A^t}) \to \H^1(k,A) \to \H^1(k,\Pic_{A^t})$. The automorphism
$-1$ of $A^t$ acts by $1$ on $\H^0(k,NS_{A^t})$ and by $-1$ on
$\H^1(k,A)$. Thus the image of the first map, which is the kernel of the
second map, is $2$-torsion. This image can, however, be non-zero -- this
phenomenon is explored in much greater depth in \cite{PoSt} where,
amongst other things, the non-triviliality of this image is identified
to be the reason (over a global field $k$) the order of $\Sha(A)$ can be
a non-square, assuming it is finite.} and the second map is injective.
Thus, if $X$ has odd order in $\H^1(k,A)$, then its image $\alpha_X$
under the map $\alpha$ has the same order. Thus, it suffices to show
that $\ind(\alpha_X) | \per(X)^g$.

Employing the notation of Lemma \ref{descent lemma}, we know that there
exists a line bundle $M$ on $X$ whose class in $\NS_{X/S} \cong
\NS_{A/S}$ is the class of $L^{\per(X)}$. In this situation, using the
notation of Theorem \ref{main thm}, I claim that $\F(M)$ is an
$\alpha_X$-twisted vector bundle, and that it has rank $\per(X)^g$. Both
these assertions are fppf local on the base and, therefore, we may
assume that $k$ is algebraically closed and, consequently, that $X$ is
trivial and $\F$ is the regular Fourier-Mukai functor, upto tensoring by
line bundles on $A^t$ (which is the same as changing the choice of
trivialisation of $X$). Since the property of being a vector bundle of a
specified rank is invariant under tensoring by line bundles, we may
assume that $\F$ is the regular Fourier-Mukai functor. Thus, we're
reduced to proving that if $M$ is line bundle on $A$ whose class in
$\NS_{A/S}$ is the same as that of $L^n$, then $\F(M)$ is a vector
bundle, and has rank $n^g$.

Since $M$ is non-degenerate, $\F(M)$ is a vector bundle (upto a
translation in the derived category). To check that it has the right
rank, recall the Fourier-Mukai functor turns Euler characteristics into
ranks, upto a sign. Thus, we need to show $\chi(M) = \pm n^g$. We know
that $\chi(M)^2 = \#K(M)$, where $K(M)$ is the kernel of the isogeny
$\phi_M:A \to A^t$. Since $\phi_M = \phi_{L^n}$, we have $\chi(M)^2 =
\#K(L^n)$. Since $L$ is a principal polarisation, $K(L)$ is trivial.
Therefore, $\#K(L^n) = A[n]$ which has rank $n^{2g}$. Thus, $\chi(M) =
\pm n^g$, as was required.
\end{proof}
\vfill
\pagebreak

\end{document}